\documentclass{amsart}%
\usepackage{amsfonts}
\usepackage[hiresbb]{graphicx}
\usepackage{amsmath}
\usepackage{amssymb}
\usepackage{version}
\usepackage{epsfig}
\usepackage{epstopdf}
\usepackage{bibmods}%
\providecommand{\U}[1]{\protect\rule{.1in}{.1in}}
\DeclareGraphicsExtensions{.eps}
\theoremstyle{theorem}
\newtheorem{defn}{Definition}[section]
\newtheorem{thm}[defn]{Theorem}
\newtheorem{exmp}[defn]{Example}
\newtheorem{lem}[defn]{Lemma}
\newtheorem{rmk}[defn]{Remark}
\newtheorem{prop}[defn]{Proposition}

\newtheorem{notat}[defn]{Notation}

\newenvironment{pf}[1][Proof]{\textbf{#1.} }{\ \rule{0.5em}{0.5em}}

\begin{document}
\title[Invariance of bifurcation equations]{Invariance of bifurcation equations for high degeneracy bifurcations of
non-autonomous periodic maps}
\author[H. Oliveira]{Henrique M. Oliveira}
\address{Center for Mathematical Analysis, Geometry and Dynamical Systems, Mathematics
Department, Instituto Superior T\'{e}cnico, Universidade de Lisboa, Av.
Rovisco Pais, 1049-001 Lisboa, Portugal }
\email{holiv@math.ist.ulisboa.pt}
\date{February, 2015 - AMS: Primary:37G15; Secondary: 39A28}

\begin{abstract}
Bifurcations of the $A_{\mu}$ class in Arnold's classification, in
non-autonomous $p$-periodic difference equations generated by parameter
depending families with $p$ maps, are studied. It is proved that the
conditions of degeneracy, non-degeneracy and unfolding are invariant relative
to cyclic order of compositions for any natural number $\mu$. The main tool
for the proofs is local topological conjugacy. Invariance results are
essential to the proper definition of the bifurcations of class $A_{\mu}$, and
lower codimension bifurcations associated, using all the possible cyclic
compositions of the fiber families of maps of the $p$-periodic difference
equation. Finally, we present two actual examples of class $A_{3}$ or
swallowtail bifurcation occurring in period two difference equations for which
the bifurcation conditions are invariant.

\end{abstract}
\maketitle

\section{Introduction}

\subsection{Motivation}

Some works on bifurcation theory of non-autonomous dynamical systems emerged
recently, \cite{ELO,POT3}. There are some difficulties to overcame in
non-autonomous systems, both with continuous or discrete time. As a starting
point it is necessary to set a proper definition of dynamical system
\cite{BHS,ES0,Kol} and of attractor and repeller \cite{AUL}. It is also
necessary to define clearly the concept of bifurcation. There is a good set of
research works on this subject such as
\cite{AlvSil,HUL,KL1,KL2,KL0,Kol,POT0,POT1,POT3,POT2,RAS0,RAS1}.

In this paper we are concerned with the definition of bifurcation equations
for local bifurcations in one-dimensional $p$-periodic maps or $p$-periodic
difference equations. In particular, we focus our attention on the $A_{\mu}$
class of bifurcations in the classification of Arnold \cite{AR0,AR} for the
positive integer $\mu$. The main result of the paper is to prove the
invariance of the $A_{\mu}$ bifurcation conditions relative to the cyclic
order of the maps in the iteration. Actually, we prove all the results for
alternating maps, i.e., with $p=2$ or two fiber maps and for fixed points of
the composition maps. This approach as the advantage of being simple in
presentation, notation and comfortable to the reader compared to the direct
study of $p$ compositions and general $k$-periodic orbits. The generalization
to periodic orbits of $p$-periodic maps is carried after, being only an
exercise of composition and repeated application of the methods developed for
alternating maps.

The $A_{\mu}$ class of bifurcation in the autonomous case occurs when one has
one real dynamic variable $x$, the parameter space is real $\mu$-dimensional
and the related family of mappings satisfy a set of degeneracy conditions.
These conditions provide topological equivalence to the unfolding of the germ
\cite{AR} $x\pm x^{\mu+1}$ at the origin. Since there are many different
approaches in the literature, in this work we follow the definitions of
\cite{AR} concerning germ, topological equivalence, unfolding, codimension and
classification of singularities and bifurcations. We suggest as an
introduction to the general subject of bifurcations the book \cite{KU}. The
$A_{\mu}$ class includes the fold, for $\mu=1$, the cusp, for $\mu=2$, the
swallowtail, for $\mu=3$, and the butterfly, for $\mu=4$ \cite{AR,Gil,TH,ZE}.

\begin{figure}[h]
\begin{center}
\includegraphics[
height=2.3375in,
width=2.5062in
]
{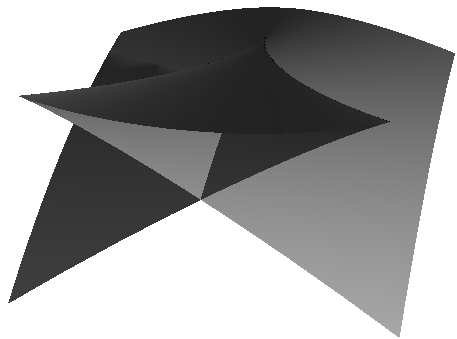}
\end{center}
\caption{The bifurcation set in the parameter space near the origin, the most
degenerate point where the $A_{3}$ singularity occurs. The control space is
real three dimensional. The cut facing the observer is at $a=1$.}%
\label{FIG9}%
\end{figure}

At the end of this paper we use the equations of the swallowtail bifurcation,
i.e., $A_{3}$ as an example of our results. In this case the bifurcation
set\footnote{For bifurcation set see definition \ref{BifSet} below.} in
parameter space is made up of three surfaces of fold bifurcations, which meet
in two lines of cusp bifurcations and one line of simultaneous double fold,
which in turn meet at a single swallowtail bifurcation point as we can see in
Figure \ref{FIG9}. This bifurcation has codimension three \cite{KU}, since one
needs three independent parameters to unfold completely the bifurcation.

On the subject of codimension see also \cite{AR,CH,GS,IOO}; we note that the
definition of codimension of \cite{GS} is different from the one provided by
\cite{AR} and \cite{KU} but the results are basically the same, modulus
personal gusto.

The $p$ maps of the family can exhibit a plethora of geometrical behavior not
present when we study lower codimension bifurcation. For instance, for $\mu=3$
and alternating maps, the Schwarzian derivative cannot be simultaneously
negative for the two maps at the singularity, as will see in the last section.
The negative Schwarzian condition restricts severely the geometry of the
families of mappings \cite{DS,SR}. Without the negative Schwarzian, we have in
the unfolding of this singularity a variety of dynamic phenomena not usually
seen in most of the works in one dimensional discrete dynamics \cite{BB,DS}.
Obviously, in this scenario one does not benefit from Singer's Theorem
\cite{S}.

\subsection{Overview}

We organized this paper in four sections. In Section \ref{SectionPrelim} we
introduce basic concepts including a brief recollection of the $A_{\mu}$
bifurcation equations for families of autonomous real maps.

In Section \ref{SectSwallow}, the core of this work, we study in detail the
equations of bifurcation for alternating systems. We prove that when we
perform a change in the order of composition of the maps the degeneracy
conditions, the non-degeneracy conditions and the transversality conditions
remain invariant. These results establish that this type of bifurcation is
well defined in the general case of alternating systems. Finally we provide a
straightforward generalization of the results on alternating system to
periodic orbits of $p$-periodic maps.

In Section \ref{Examples}\ we prove some conditions that override the
possibility of $A_{3}$ or swallowtail bifurcation in the case of alternating
maps, where $p=2$. Finally, we present two examples concerning alternating
maps. These examples show that this class of high degeneracy bifurcations
occurs in simple applications without the need of exotic constructions.

\section{\label{SectionPrelim}Preliminaries}

\subsection{Basic definitions and notation}

We define non-autonomous dynamical system using the definitions of \cite{Kol}.
Consider the non-autonomous iteration given by%
\begin{equation}
x_{n+1}=f_{n}(x_{n})\text{, }x_{n}\in I_{n}\text{, with }n\in\mathbb{%
%TCIMACRO{\U{2115} }%
%BeginExpansion
\mathbb{N}
%EndExpansion
}\text{,} \label{pp}%
\end{equation}
where $I_{n}$ are real intervals (not necessarily compact and maybe $%
%TCIMACRO{\U{211d} }%
%BeginExpansion
\mathbb{R}
%EndExpansion
$) which are the fibers of the non-autonomous dynamical system. The usual
distance is defined in $I_{n}$. The iteration starts at the initial condition
$x_{0}\in I_{0}$. Each map $f_{n}$ is at least continuous and defined such
that%
\[%
\begin{array}
[c]{cccc}%
f_{n}: & I_{n} & \longrightarrow & I_{n+1}\text{,}\\
& x_{n} & \longmapsto & f_{n}\left(  x_{n}\right)
\end{array}
\]
and $f_{n}\left(  I_{n}\right)  \subseteq I_{n+1}$.

The system is periodic of period $p$ if
\[
f_{n+p}=f_{p}\text{ and }f_{n+p}\left(  I_{n}\right)  \subseteq I_{n}\text{,}%
\]
for every $n\in%
%TCIMACRO{\U{2115} }%
%BeginExpansion
\mathbb{N}
%EndExpansion
$, where $p$ is a minimal positive integer. When $p=2$ with the fibers $I_{0}$
and $I_{1}$, such that%
\[
f_{0}\left(  I_{0}\right)  \subseteq I_{1}\text{ and }f_{1}\left(
I_{1}\right)  \subseteq I_{0}\text{,}%
\]
we say that we have an \textit{alternating system}.

In the sequel, by ${\mathcal{C}}\left(  I\right)  $ we denote the collection
of all continuous maps in its domain $I$, by ${\mathcal{C}}^{1}\left(
I\right)  $ the collection of all continuously differentiable elements of
${\mathcal{C}}\left(  I\right)  $ and, in general by ${\mathcal{C}}^{s}\left(
I\right)  ,$ $s\geq1,$ the collection of all elements of ${\mathcal{C}}\left(
I\right)  $ having continuous derivatives up to order $s$ in $I.$

Let $f\in{\mathcal{C}}^{1}\left(  I\right)  $, and let $q$ be a periodic point
of period $m$. Denoting the derivative by $D$, $q$ is called a
\emph{hyperbolic attractor} if $|Df^{m}(q)|<1$, a \emph{hyperbolic repeller}
if $|Df^{m}(q)|>1$, and \emph{non-hyperbolic} if $|Df^{m}(q)|=1$.

\begin{defn}
We say that two continuous maps $f:I\rightarrow I$ and $g:J\rightarrow J$, are
topologically conjugate, if there exists a homeomorphism $h:I\rightarrow J$,
such that $h\circ f=g\circ h$. We call $h$ the topological conjugacy of $f$
and $g$.
\end{defn}

We use the capital $\Lambda$ for a vector parameter in $\mathbb{R}^{\mu}$.

\begin{defn}
\label{BifSet}If $f_{\Lambda}$ is a family of maps, then the regular values
$\Lambda$ of the parameters are those which have the property that
$f_{\widetilde{\Lambda}}$ is topologically conjugate to $f_{\Lambda}$ for all
$\widetilde{\Lambda}$ in some open neighbourhood of $\Lambda$. If $\Lambda$ is
not a regular value, it is a \emph{bifurcation value}. The collection of all
the bifurcation values is the \emph{bifurcation set}, $\Omega\subset
\mathbb{R}^{\mu}$, in the parameter space.
\end{defn}

Let $f_{\Lambda}$ be a family of maps in ${\mathcal{C}}^{s}\left(  I\right)
$. Let $\Lambda_{0}$ be a particular vector parameter and $a\in I$ be a fixed
point of $f_{\Lambda_{0}}$, i.e.,%
\[
a=f_{\Lambda_{0}}\left(  a\right)  \text{,}%
\]
the condition of $a$ being non-hyperbolic is necessary for the existence of a
local bifurcation. The existence and nature of that bifurcation depends on
other symmetry and differentiable conditions that we will see bellow. If there
exists a local bifurcation we say that $(a,\Lambda_{0})$ is a
\emph{bifurcation point }(when there is no risk of confusion, we say that $a$
is a \emph{bifurcation point}).

\begin{notat}
For notational simplicity we consider the real vector parameter $\Lambda$ as a
standard variable along with the dynamic variable $x$, i.e., we write
\[
f_{\Lambda}\left(  x\right)  =f\left(  x,\Lambda\right)  ,
\]
keeping in mind that the compositions are always in the dynamic variable $x$.

When there are no danger of confusion and no operations regarding the
parameter, we denote the evaluation of functions depending on the dynamic
variable and the parameter omitting the later, for instance $f_{\Lambda
}\left(  x\right)  =f\left(  x,\Lambda\right)  $ will be denoted by $f\left(
x\right)  $ in order to avoid overload the complicated notation needed for the
computations of high order chain rules. Nevertheless, all the maps in this
paper depend on the parameter as well on the dynamic variable. We deal with
parameter depending families of maps, even when that dependence is not visible
in some formulas or expressions.

We denote the derivatives relative to some variable $y$ by $\partial_{y}$.
Repeated differentiation relative to the same variable is denoted by
$\partial_{y^{n}}$, for instance $\partial_{yyy}=\partial_{y^{3}}$. When there
is no danger of confusion, we denote strict partial derivatives, i.e., not
seeing composed functions, by a subscript. For instance, the third partial
derivative of $f$ relative to $y$ is, in that case, denoted by $f_{yyy}$ or
$f_{y^{3}}$.
\end{notat}

This means, in particular, that when dealing with the composition of real
scalar functions $g\left(  x,t\right)  $ and $f\left(  x,t\right)  $, such
that $g\circ f\left(  x,t\right)  =g\left(  f\left(  x,t\right)  ,t\right)  $,
we have the chain rules%
\begin{align*}
\partial_{t}g\left(  f\left(  x,t\right)  ,t\right)   &  =g_{x}\left(
f\left(  x,t\right)  ,t\right)  f_{t}\left(  x,t\right)  +g_{t}\left(
f\left(  x,t\right)  ,t\right)  \text{,}\\
\partial_{x}g\left(  f\left(  x,t\right)  ,t\right)   &  =g_{x}\left(
f\left(  x,t\right)  ,t\right)  f_{x}\left(  x,t\right)  \text{.}%
\end{align*}

Along this paper we deal with $p$-periodic sequences of maps $f_{0}$, $f_{1}%
$,\ldots\ , $f_{p-1}$ on a real dynamic variable $x$ and depending on a real
vector parameter $\Lambda$, such that%
\[%
\begin{array}
[c]{cccc}%
f_{j}: & I_{j}\times\Theta & \longrightarrow & I_{j+1}\\
& \left(  x,\Lambda\right)  & \longmapsto & f_{0}\left(  x,\Lambda\right)
\end{array}
\]
for $j=0,\ldots,p-1$. The fibers $I_{j}$ for the dynamic variable are
intervals of $%
%TCIMACRO{\U{211d} }%
%BeginExpansion
\mathbb{R}
%EndExpansion
$ and $\Theta\subset%
%TCIMACRO{\U{211d} }%
%BeginExpansion
\mathbb{R}
%EndExpansion
^{\mu}$ is the parameter set, $f_{j}\in\mathcal{C}^{\mu+1}\left(
I_{j}\right)  $ and $f_{j}\in\mathcal{C}^{1}\left(  \Theta\right)  $, with
$\mu$ a positive integer. Moreover, the property%
\[
f_{j}\left(  I_{j},\Lambda\right)  \subseteq I_{j+1\left(  \operatorname{mod}%
p\right)  }\text{, holds for all }\Lambda\in\Theta\text{.}%
\]
In this paper we use the convention that capital letters are used for
compositions of maps in the dynamic variable. Capital $F$ and $G$ will be
always used for direct and reverse composition of alternating maps%
\[
F=f_{1}\circ f_{0}\text{ and }G=f_{0}\circ f_{1}.
\]
Consider the set of indexes $j=0,\ldots,p-1$ for the $p$-periodic system. We
set the following notation for the $p$ compositions
\begin{align*}
F_{0}  &  =f_{p-1}\circ\cdots\circ f_{0},\\
F_{1}  &  =f_{0}\circ f_{p-1}\circ\cdots\circ f_{1},\\
&  \vdots\\
F_{p-1}  &  =f_{p-2}\circ\cdots f_{0}\circ f_{p-1}\text{.}%
\end{align*}
Repeated composition (always in the dynamic variable) is denoted by%
\[
f^{k}=\underset{k}{\underbrace{(f\circ\cdots\circ f)}}\text{,}%
\]
where $k$ is a positive integer.

\subsection{Conditions for the $A_{\mu}$ class of bifurcations in autonomous
systems}

In this paragraph, we recall briefly the conditions of class $A_{\mu}$ of
local bifurcations in Arnold classification as explained in Theorem of page 20
in Arnold et al. \cite{AR}. For the iteration of maps, the normalized germ of
the class $A_{\mu}$ is $x\pm x^{\mu+1}$ and has the principal family
\cite{AR}, also called prototype polynomial or normal form \cite{KU}%
\[
x\pm x^{\mu+1}+\lambda_{1}+\ldots+\lambda_{\mu}x^{\mu-1}\text{,}%
\]
where $\lambda_{j}$, $j=1,\ldots,\mu$, are real parameters.

Giving an autonomous discrete dynamical system generated by the iteration of
$f$, in order to compute the bifurcation points of class $A_{\mu}$, one has to
solve the bifurcation equations \cite{KU}%
\begin{equation}%
\begin{array}
[c]{l}%
{f}(x,\Lambda)=x\text{, fixed point equation}\\
{f}_{x}(x,\Lambda)=1\text{, non-hyperbolicity condition.}%
\end{array}
\label{Fold}%
\end{equation}
The simplest of such local bifurcations is the saddle node bifurcation, i.e.,
$A_{1}$. One assumes, in this case, the generic non-degeneracy condition
\begin{equation}
{f_{xx}}(x,\lambda)\not =0 \label{Nondeg}%
\end{equation}
and the transversality condition \cite{KU}%
\begin{equation}
f_{\lambda}\left(  x,\lambda\right)  \not =0,\text{ with }\lambda\in%
%TCIMACRO{\U{211d} }%
%BeginExpansion
\mathbb{R}
%EndExpansion
\text{.} \label{unfold1}%
\end{equation}
We set generically that $\lambda\in%
%TCIMACRO{\U{211d} }%
%BeginExpansion
\mathbb{R}
%EndExpansion
$, since one needs only one parameter to unfold locally this singularity
\cite{A,AR,CH,GS,GU,KU}. The normalized germ of this bifurcation is%
\[
x\pm x^{2}\text{,}%
\]
with principal family%
\[
x\pm x^{2}+\lambda\text{,}%
\]
which is weak topologically conjugated to any other family \cite{AR,KU}
satisfying the bifurcation conditions.

Adding degeneracy conditions, one obtains higher degeneracy local bifurcations.

Therefore, the equations for the occurrence of $A_{\mu}$ class of bifurcations
for a general positive integer $\mu$ are%
\begin{equation}%
\begin{array}
[c]{l}%
{f}(x,\Lambda)=x,\\
{f_{x}}(x,\Lambda)=1,\\
{{f}_{xx}}(x,\Lambda)=0,\\
\vdots\\
{f_{x^{\mu}}}(x,\Lambda)=0,
\end{array}
\label{DEGswallow}%
\end{equation}
with solution $\left(  a,\Lambda_{0}\right)  $. It is easy to see that these
conditions are satisfied by the normalized germ $x\pm x^{\mu+1}$ at the
origin. One has the non-degeneracy condition%
\begin{equation}
{f_{x^{\mu+1}}}(a,\Lambda_{0})\not =0, \label{NDEGswallow}%
\end{equation}
for $\mu=2$ we have the cusp, for $\mu=3$ the swallowtail and for $\mu=4$ the
butterfly \cite{AR,CH,EH,GS,GU,KU}. The transversality condition (see pages
66, 297, 298 and 303 of \cite{KU}) at the solution of the above conditions is
given by the condition on the non-singularity of the Jacobian matrix of the
map $\left(  f,f_{x},f_{xx},\ldots,f_{x^{\mu-1}}\right)  $ relative to the
parameters at the bifurcation point%
\begin{equation}
\det\left[
\begin{array}
[c]{lll}%
f_{\lambda_{1}}(a,\Lambda_{0}) & \cdots & f_{\lambda_{\mu}}(a,\Lambda_{0})\\
\vdots & \ddots & \vdots\\
f_{x^{\mu-1}\lambda_{1}}(a,\Lambda_{0}) & \cdots & f_{x^{\mu-1}\lambda_{\mu}%
}(a,\Lambda_{0})
\end{array}
\right]  \not =0, \label{UnfoldSwallow}%
\end{equation}
and assures that the vector parameter is enough to unfold the local
bifurcation \cite{KU}. This happens since condition (\ref{UnfoldSwallow})
assures that the $\mu$ lower order terms in the Taylor polynomial of $f$
depend uniquely on the $\mu$ components of $\Lambda$, i.e., $\lambda_{1}%
$,\ldots\ ,$\lambda_{\mu}$.

\section{\label{SectSwallow}$A_{\mu}$ class of bifurcation in families of
$p$-periodic maps}

\subsection{Invariance of the bifurcation conditions}

\subsubsection{On the invariance of the degeneracy and non-degeneracy
conditions for alternating systems\label{Invariance}}

In this paragraph, we study the invariance of the degeneracy conditions of
alternating families of maps for all the singularities of class $A_{\mu}$,
using topological conjugacy.

Given an initial condition $x_{0}\in I_{0}$ the alternating system is given by
the iteration%
\begin{equation}
x_{n+1}=f_{n\left(  \operatorname{mod}2\right)  }\left(  x_{n},\Lambda\right)
\text{, }x_{n}\in I_{n\left(  \operatorname{mod}2\right)  }. \label{it}%
\end{equation}

If there is a pair $\left(  a,I_{0}\right)  ,$ such that after $2$ iterations,
the iteration returns to $\left(  a,I_{0}\right)  $ we say that $a$ is a
\emph{periodic point in the fiber }$I_{0}$ with period $2$. We note that the
point $b=f_{0}\left(  a,\Lambda\right)  \ $is also a periodic point in the
fiber $I_{1}$ with period $2$. Consider the compositions $F$ and $G$, we have%
\[
a=F\left(  a,\Lambda\right)  \text{ and }b=G\left(  b,\Lambda\right)  \text{.}%
\]
In other words: $a$ (resp. $b$) is a periodic point with period $2$ in fiber
$I_{0}$ (resp. $I_{1}$) of the alternating system\ (\ref{it})\ iff $a$ (resp.
$b$) is a fixed point of $F$ (resp. $G$).

These below are the bifurcation equations with $\mu-1$ degeneracy conditions
on derivatives on $x$ stated for $F$ and $G$ which are exactly the same as in
the non-autonomous case
\begin{equation}
\left\{
\begin{array}
[c]{l}%
{F}(x,\Lambda)=x,\\
{F_{x}}(x,\Lambda)=1,\\
{{F}_{xx}}(x,\Lambda)=0,\\
\vdots\\
{F_{x^{\mu}}}(x,\Lambda)=0,
\end{array}
\right.  \text{ and }\left\{
\begin{array}
[c]{l}%
{G}(x,\Lambda)=x,\\
{G_{x}}(x,\Lambda)=1,\\
{{G}_{xx}}(x,\Lambda)=0,\\
\vdots\\
{G_{x^{\mu}}}(x,\Lambda)=0.
\end{array}
\right.  \label{bif}%
\end{equation}
These equations have different solutions for the dynamic variable $x$,
depending on the fiber we choose. At the solutions of (\ref{bif}), the
non-degeneracy conditions are
\begin{equation}
{F_{x^{\mu+1}}}(a,\Lambda_{0})\not =0\text{ and }{G_{x^{\mu+1}}}(b,\Lambda
_{0})\not =0 \label{nond}%
\end{equation}

A natural question arises:

\textbf{Are the solutions in the parameter space equal for the different
compositions }$F$ and $G$\textbf{?}

The same query was posed in \cite{EH,ELO} and positively solved in the
particular cases dealt in those works for degeneracy conditions until the
cusp, i.e., $\mu=2$.

Indeed, in this paragraph we show that the answer to the question is positive
in the general case. We prove that if a parameter vector satisfies equations
(\ref{bif}) for $F$ then it is a solution of the system for $G$.

The next lemma will be used to solve the general problem of the symmetry of
the bifurcation equations with respect to the order of composition.

\begin{lem}
\label{lema}Let $\mu\geq1$ and let $h$ and $f$ be real functions satisfying
the conditions:

\begin{enumerate}
\item \label{A1} there exists $a$ such that $f(a)=a$ and $f$ is a Lipschitz
homeomorphism in some open interval $I$ containing $a$;

\item \label{A2}$h$ is a Lipschitz homeomorphism with Lipschitz constant $L$
in a open neighborhood $I_{h}$ of $a$, and there exists an open neighborhood,
$J_{b}$, of $h(a)=b$ such that its inverse, $h^{-1}$, is also Lipschitz
continuous with Lipschitz constant $M$;

\item \label{A3}$\ $%
\[
\underset{x\rightarrow a}{\lim}\frac{\left\vert f\left(  x\right)
-x\right\vert }{\left\vert x-a\right\vert ^{\mu}}=0\text{ and }\underset
{x\rightarrow a}{\lim}\frac{\left\vert f\left(  x\right)  -x\right\vert
}{\left\vert x-a\right\vert ^{\mu+1}}>0
\]

\end{enumerate}

Then $g$, the conjugate of $f$ by the homeomorphism $h$,%
\[
g=h\circ f\circ h^{-1}%
\]
satisfies%
\begin{equation}
\underset{y\rightarrow b}{\lim}\frac{\left\vert g\left(  y\right)
-y\right\vert }{\left\vert y-b\right\vert ^{\mu}}=0\text{ and }\underset
{y\rightarrow b}{\lim}\frac{\left\vert g\left(  y\right)  -y\right\vert
}{\left\vert y-b\right\vert ^{\mu+1}}>0 \label{lim}%
\end{equation}

\end{lem}

\begin{pf}
We first compute the domain $J$ where $y$ ranges when we compute the limit
(\ref{lim}). Of course, we take $J$ to be an open interval containing $b$ such
that $J\subseteq h\left(  I_{h}\right)  $. The limit has meaning if $J$ also
satisfies
\begin{align}
h^{-1}\left(  J\right)   &  \subseteq I\label{set1}\\
f\left(  h^{-1}\left(  J\right)  \right)   &  \subseteq I_{h}.\nonumber
\end{align}
As both $f$ and $h^{-1}$ are homeomorphisms we can choose the open interval
$J$ small enough just to satisfy the conditions (\ref{set1}). We note that
$a\in h^{-1}\left(  J\right)  $ and $a\in f\left(  h^{-1}\left(  J\right)
\right)  $.

Let us consider the limit (\ref{lim}).\ We have%
\begin{align*}
0\leq\underset{y\rightarrow b}{\lim}\frac{\left\vert g\left(  y\right)
-y\right\vert }{\left\vert y-b\right\vert ^{\mu}}  &  =\underset{y\rightarrow
b}{\lim}\frac{\left\vert (h\circ f\circ h^{-1})\left(  y\right)  -(h\circ
h^{-1})\left(  y\right)  \right\vert }{\left\vert y-b\right\vert ^{\mu}}\\
&  \leq L\underset{y\rightarrow b}{\lim}\frac{\left\vert (f\circ
h^{-1})\left(  y\right)  -h^{-1}\left(  y\right)  \right\vert }{\left\vert
y-b\right\vert ^{\mu}}\\
&  =L\underset{y\rightarrow b}{\lim}\frac{\left\vert (f\circ h^{-1})\left(
y\right)  -h^{-1}\left(  y\right)  \right\vert }{\left\vert h^{-1}\left(
y\right)  -a\right\vert ^{\mu}}\left(  \frac{\left\vert h^{-1}\left(
y\right)  -a\right\vert }{\left\vert y-b\right\vert }\right)  ^{\mu}\\
&  =L\underset{y\rightarrow b}{\lim}\frac{\left\vert (f\circ h^{-1})\left(
y\right)  -h^{-1}\left(  y\right)  \right\vert }{\left\vert h^{-1}\left(
y\right)  -a\right\vert ^{\mu}}\left(  \frac{\left\vert h^{-1}\left(
y\right)  -h^{-1}\left(  b\right)  \right\vert }{\left\vert y-b\right\vert
}\right)  ^{\mu}\\
&  \leq LM^{\mu}\underset{y\rightarrow b}{\lim}\frac{\left\vert (f\circ
h^{-1})\left(  y\right)  -h^{-1}\left(  y\right)  \right\vert }{\left\vert
h^{-1}\left(  y\right)  -a\right\vert ^{\mu}}.
\end{align*}
As
\[
\underset{y\rightarrow b}{\lim}h^{-1}\left(  y\right)  =a,
\]
if we set $h^{-1}(y)=x$, it follows
\[
0\leq\underset{y\rightarrow b}{\lim}\frac{\left\vert g\left(  y\right)
-y\right\vert }{\left\vert y-b\right\vert ^{\mu}}\leq LM^{\mu}\underset
{x\rightarrow a}{\lim}\frac{\left\vert f\left(  x\right)  -x\right\vert
}{\left\vert x-a\right\vert ^{\mu}}=0.
\]
When one has $\underset{x\rightarrow a}{\lim}\frac{\left\vert f\left(
x\right)  -x\right\vert }{\left\vert x-a\right\vert ^{\mu+1}}>0$, we apply
similar reasoning to $f$ to get%
\[
0<\underset{x\rightarrow a}{\lim}\frac{\left\vert f\left(  x\right)
-x\right\vert }{\left\vert x-a\right\vert ^{\mu+1}}\leq ML^{\mu+1}%
\underset{y\rightarrow b}{\lim}\frac{\left\vert g\left(  y\right)
-y\right\vert }{\left\vert y-b\right\vert ^{\mu+1}},
\]
therefore%
\[
\underset{y\rightarrow b}{\lim}\frac{\left\vert g\left(  y\right)
-y\right\vert }{\left\vert y-b\right\vert ^{\mu+1}}>0,
\]
as desired.
\end{pf}

Using the previous lemma we can easily prove the next result.

\begin{thm}
\label{t1} Let $\mu\geq2$ and let be the alternating family of maps with
individual mappings $f_{0}$, $f_{1}$ with $f_{0}\in\mathcal{C}^{\mu+1}\left(
I_{0}\right)  $, $f_{1}\in\mathcal{C}^{\mu+1}\left(  I_{1}\right)  $ in the
dynamic variable and the compositions $F=f_{1}\circ f_{0}$ and $G=f_{0}\circ
f_{1}$, satisfying:

\begin{enumerate}
\item \label{D1}There exist $a$, $b$, fixed points of $F\ $and $G$
respectively
\begin{align*}
a  &  =F\left(  a,\Lambda\right)  \text{,}\\
b  &  =G\left(  b,\Lambda\right)  \text{.}%
\end{align*}

\item \label{D2}The non-hyperbolicity condition for $F$%
\[
\left.  \partial_{x}F\left(  x,\Lambda\right)  \right\vert _{x=a}=\left.
\partial_{x}f_{1}\left(  x,\Lambda\right)  \right\vert _{x=b}\left.
\partial_{x}f_{0}\left(  x,\Lambda\right)  \right\vert _{x=a}=1.
\]

\item \label{D3}Higher degeneracy conditions for $F$%
\[
\ \left.  \partial_{x^{i}}F\left(  x,\Lambda\right)  \right\vert
_{x=a}=0,\mathit{\ }\text{\textit{for every }}2\leq i\leq\mu.
\]

\item \label{D4}The non-degeneracy condition for $F$%
\[
\ \left.  \partial_{x^{\mu+1}}F\left(  x,\Lambda\right)  \right\vert
_{x=a}\not =0.
\]

\end{enumerate}

Then, the composition $G$, satisfies%
\[
\left.  \partial_{x^{i}}G\left(  x,\Lambda\right)  \right\vert _{x=b}=0,\text{
for every }2\leq i\leq\mu
\]
and
\[
\left.  \partial_{x^{\mu+1}}G\left(  x,\Lambda\right)  \right\vert
_{x=b}\not =0.
\]

\end{thm}

\begin{pf}
Properties (\ref{D1}) and (\ref{D2}) imply that $f_{0}$, $f_{1}$ are
diffeomorphisms in suitable neighborhoods of $a$ and $b$, respectively.
Therefore, we can define local inverses. Being local diffeomorphisms, $f_{0}$
and $f_{1}$ are also local Lipschitz continuous and so their inverses. In
particular $f_{0}\left(  a\right)  =b$ and $f_{0}^{-1}\left(  b\right)  =a$.
We apply Lemma \ref{lema} to $F$ and $G$ making the identification $f=F$,
$g=G$ and $h=f_{0}$.

By (\ref{D2}) and (\ref{D3})%
\[
\lim_{x\rightarrow a}\frac{\left\vert F\left(  x\right)  -x\right\vert
}{\left\vert x-a\right\vert ^{\mu}}=\lim_{x\rightarrow a}\frac{\left\vert
\left(  F\left(  x\right)  -x\right)  -\left(  F\left(  a\right)  +a\right)
\right\vert }{\left\vert x-a\right\vert ^{\mu}}=0.
\]

Therefore, $h=f_{0}$ and $F$ satisfy the hypotheses of Lemma \ref{lema}, and
hence the thesis with $G\left(  x\right)  =(f_{0}\circ F\circ f_{0}%
^{-1})\left(  x\right)  $. Thus, we obtain%
\[
\lim_{x\rightarrow b}\frac{\left\vert \left(  f_{0}\circ F\circ f_{0}%
^{-1}\right)  \left(  x\right)  -x\right\vert }{\left\vert x-b\right\vert
^{\mu}}=\lim_{x\rightarrow b}\frac{\left\vert G\left(  x\right)  -x\right\vert
}{\left\vert x-b\right\vert ^{\mu}}=0
\]
and%
\[
\lim_{x\rightarrow b}\frac{\left\vert \left(  f_{0}\circ F\circ f_{0}%
^{-1}\right)  \left(  x\right)  -x\right\vert }{\left\vert x-b\right\vert
^{\mu+1}}=\lim_{x\rightarrow b}\frac{\left\vert G\left(  x\right)
-x\right\vert }{\left\vert x-b\right\vert ^{\mu+1}}>0,
\]
that is, the first $\mu$ derivatives of $G\left(  x\right)  -x$ are zero at
$b$ and the non-degeneracy condition holds as well.
\end{pf}

%As in \cite{CL} and in \cite{EH}, given two maps $f_0  \in \mathcal{C} (I_{0})$
%and $f_1 \in\mathcal{C} (I_{1})$, we denote by $[f_{1}, f_{2}]$
%the sequences $\{x_n\}$ with $x$ ranging in $I_{1}$, defined by the rule
%$\{x,f_{1}(x), f_{2}(f_{1}(x)), f_{1}(f_{2}(f_{1}(x))), \ldots\}$, i.e. we set
%$x_{0}=x$ and $x_{n+1}= f_{i}(x_{n})$, for every $n \geq 0$, where $i=1$ if $n$ is
%even and $i=2$ if $n$ is odd, that is the {\it the non autonomous system generated by
%$f_1$ and $f_2$} or {\it the alternating dynamics generated by $f_1$ and $f_2$}. \\

\subsubsection{Example using Fa\`{a} di Bruno Formula}

Although seeming needless after the previous results, the next example will be
important to deduce properties on the geometrical behavior of the composition
of maps related to the swallowtail bifurcation at the beginning of Section
\ref{Examples}. On the other hand, it is interesting to recover the Fa\`{a} di
Bruno's formula \cite{JO,LA}, since we will use it to prove the invariance of
the transversality conditions. We think that it is possible to establish
combinatorial results, using Theorem \ref{t1}\ on both ends of the general
formula for the derivatives of the compositions. This is an interesting open
line of research for readers interested in Bell polynomials and other relevant
combinatorial quantities associated with the Fa\`{a} di Bruno's Formula, see
\cite{JO,NO,RO}.

\begin{exmp}
\label{example}(Alternating maps) Giving two real maps $f$ and $g$ defined in
real intervals $I_{0}$ and $I_{1}$ we prove directly that if the second
derivative relative to the dynamic variable of any of the two maps $g\circ f$
and $f\circ g$ is zero, then also the other must be zero, disregarding the
order of composition. The same holds for the third derivatives. We do this
directly, using the chain rule for computing the derivatives of composed maps
and its generalization, the \textit{Fa\`{a} di Bruno's Formula} \cite{JO}.

Let $f$ and $g$ be $\mathcal{C}^{3}$ functions satisfying the conditions:

\begin{enumerate}
\item \label{G1}$(g\circ f)\left(  a\right)  =a$ and $(f\circ g)\left(
b\right)  =b$, which is $f\left(  a\right)  =b$ and $g\left(  b\right)  =a$.

\item \label{G2} $\left.  \frac{d(g\circ f)}{dx}\left(  x\right)  \right\vert
_{x=a}=g^{\prime}\left(  b\right)  f^{\prime}\left(  a\right)  =1$.

\item \label{G3}$\left.  \frac{d^{m}(g\circ f)}{dx^{m}}\left(  x\right)
\right\vert _{x=a}=0$ for $m=2,3$.
\end{enumerate}

Let us recall the formula of Fa\`{a} di Bruno for the derivatives of the
composition%
\begin{equation}
\frac{d^{m}\left(  g\circ f\right)  }{dx^{m}}\left(  x\right)  =m!\sum
g^{\left(  n\right)  }\left(  f\left(  x\right)  \right)  \overset
{m}{\underset{j=1}{%
%TCIMACRO{\dprod }%
%BeginExpansion
{\displaystyle\prod}
%EndExpansion
}}\frac{1}{\beta_{j}!}\left(  \frac{f^{\left(  j\right)  }\left(  x\right)
}{j!}\right)  ^{\beta_{j}}\text{,} \label{bruno}%
\end{equation}
where the sum is over all different solutions $\beta_{j}$ in nonnegative
integers, $\beta_{1}$, $\ldots$, $\beta_{m}$, of the linear Diophantine
equations%
\[
\underset{j=1}{\overset{m}{\sum}}j\beta_{j}=m,\text{ and }n:=\underset
{j=1}{\overset{m}{\sum}}\beta_{j}.
\]
To avoid to overload this example with indexes we use the notation used in
\cite{RO}
\[
f_{0}=f\left(  a\right)  ,\qquad f_{1}=\left.  f^{\prime}\left(  x\right)
\right\vert _{x=a},\ldots,\qquad f_{m}=\left.  f^{\left(  m\right)  }\left(
x\right)  \right\vert _{x=a},
\]%
\[
g_{0}=g\left(  b\right)  ,\qquad g_{1}=\left.  g^{\prime}\left(  x\right)
\right\vert _{x=b},\ldots,\qquad g_{m}=\left.  g^{\left(  m\right)  }\left(
x\right)  \right\vert _{x=b},
\]%
\[
\left.  \frac{d^{m}(g\circ f)}{dx^{m}}\left(  x\right)  \right\vert
_{x=a}=\left(  gf\right)  _{m},\qquad\left.  \frac{d^{m}(f\circ g)}{dx^{m}%
}\left(  x\right)  \right\vert _{x=b}=\left(  fg\right)  _{m}.
\]
With this notation, and taking into account the hypotheses \ref{G1}, \ref{G2}
and \ref{G3}, Fa\`{a} di Bruno's Formula gives
\begin{equation}
\left(  gf\right)  _{m}=m!\sum g_{n}\overset{m}{\underset{j=1}{%
%TCIMACRO{\dprod }%
%BeginExpansion
{\displaystyle\prod}
%EndExpansion
}}\frac{1}{\beta_{j}!}\left(  \frac{f_{j}}{j!}\right)  ^{\beta_{j}}
\label{bruno1}%
\end{equation}
and%
\begin{equation}
\left(  fg\right)  _{m}=m!\sum f_{n}\overset{m}{\underset{j=1}{%
%TCIMACRO{\dprod }%
%BeginExpansion
{\displaystyle\prod}
%EndExpansion
}}\frac{1}{\beta_{j}!}\left(  \frac{g_{j}}{j!}\right)  ^{\beta_{j}}
\label{bruno2}%
\end{equation}
Condition (\ref{G2}) in this notation is now%
\begin{equation}
f_{1}g_{1}=1. \label{cc2}%
\end{equation}
Let us consider the first two cases: $m=2$ and $m=3$, cusp and swallowtail.

Let $m=2$. We shall use the formula (\ref{bruno}), therefore we have to solve
the equation%
\[
\beta_{1}+2\beta_{2}=2,
\]
for all possible values of the vector $\left(  \beta_{1},\beta_{2}\right)  $
in $\mathbb{N}\times\mathbb{N}$. The only solutions are $(\beta_{1},\beta
_{2})=(0,1)$, which gives $n=1$ and $(\beta_{1},\beta_{2})=(2,0)$, which gives
$n=2$. So we have%
\begin{align}
\left(  gf\right)  _{2}  &  =2!\left(  g_{1}\frac{1}{0!}\left(  \frac{f_{1}%
}{1!}\right)  ^{0}\frac{1}{1!}\left(  \frac{f_{2}}{2!}\right)  ^{1}+g_{2}%
\frac{1}{2!}\left(  \frac{f_{1}}{1!}\right)  ^{2}\frac{1}{0!}\left(
\frac{f_{2}}{2!}\right)  ^{0}\right)  =0\label{d21}\\
&  =g_{1}f_{2}+g_{2}f_{1}^{2}=g_{1}f_{2}+\frac{g_{2}}{g_{1}^{2}}=0\nonumber
\end{align}
and%
\begin{align}
\left(  fg\right)  _{2}  &  =2!\left(  f_{1}\frac{1}{0!}\left(  \frac{g_{1}%
}{1!}\right)  ^{0}\frac{1}{1!}\left(  \frac{g_{2}}{2!}\right)  ^{1}+f_{2}%
\frac{1}{2!}\left(  \frac{g_{1}}{1!}\right)  ^{2}\frac{1}{0!}\left(
\frac{g_{2}}{2!}\right)  ^{0}\right) \label{d22}\\
&  =f_{1}g_{2}+f_{2}g_{1}^{2}=\frac{g_{2}}{g_{1}}+f_{2}g_{1}^{2}.\nonumber
\end{align}
We solve the system with equations (\ref{cc2}) and (\ref{d21}) for $g_{2}$, to
obtain
\begin{equation}
g_{2}=g_{2}(f_{1},f_{2})=-\frac{f_{2}}{f_{1}^{3}}. \label{g2}%
\end{equation}
By substituting $g_{1}=\frac{1}{f_{1}}$ and $g_{2}\left(  f_{1},f_{2}\right)
$ in (\ref{d22}), we get%
\[
\left(  fg\right)  _{2}=-f_{1}\frac{f_{2}}{f_{1}^{3}}+f_{2}\frac{1}{f_{1}^{2}%
}=0.
\]
Let $m=3$.

By Fa\`{a} di Bruno Formula, and taking into account the hypotheses, we
obtain
\begin{align}
\left(  gf\right)  _{3}  &  =g_{1}f_{3}+3g_{2}f_{1}f_{2}+g_{3}f_{1}%
^{3}\label{d31}\\
&  =g_{1}f_{3}+\frac{3g_{2}f_{2}}{g_{1}}+\frac{g_{3}}{g_{1}^{3}}=0.\nonumber
\end{align}
We solve the system with equations (\ref{cc2}), (\ref{d21}) and (\ref{d31})
for $g_{3}$, keeping in mind that $g_{1}$ and $g_{2}$ have been computed
before. Hence, we get%
\begin{equation}
g_{3}=-\frac{f_{3}}{f_{1}^{4}}-\frac{3\left(  -\frac{f_{2}}{f_{1}^{3}}\right)
f_{2}}{f_{1}^{2}}=\frac{3f_{2}^{2}}{f_{1}^{5}}-\frac{f_{3}}{f_{1}^{4}}.
\label{g3}%
\end{equation}
By replacing $g_{1}$, $g_{2}$ and $g_{3}$ by the solutions previously
obtained, we find
\begin{align}
\left(  fg\right)  _{3}  &  =f_{1}g_{3}+3f_{2}g_{1}g_{2}+f_{3}g_{1}%
^{3}\label{d32}\\
&  =f_{1}\left(  \frac{3f_{2}^{2}}{f_{1}^{5}}-\frac{f_{3}}{f_{1}^{4}}\right)
+3f_{2}\frac{1}{f_{1}}\left(  -\frac{f_{2}}{f_{1}^{3}}\right)  +\frac{f_{3}%
}{f_{1}^{3}}=0.\nonumber
\end{align}

\end{exmp}

\subsubsection{On the invariance of the degeneracy and non-degeneracy
conditions for periodic orbits of $p$-periodic systems}

What we have just shown in Paragraph \ref{Invariance}\ is the invariance of
bifurcation equations with respect to interchange in the composition of the
alternating maps. In this paragraph we generalize the results for alternating
maps to general $p$-periodic non-autonomous systems.

\begin{thm}
\label{tt1}Let $\mu\geq2$ and let be the $p$-periodic family of maps with
individual mappings $f_{0}$, $f_{1}$,\ldots\ , $f_{p-1}$ with $f_{j}%
\in\mathcal{C}^{\mu+1}\left(  I_{j}\right)  $, with one fixed $j\in\left\{
0,\ldots,p-1\right\}  $, and a periodic point $a_{j}\in I_{j}$, with period
$p$, i.e., a fixed point of $F_{j}$, satisfying:

\begin{enumerate}
\item \label{DD1}There exist $a_{0}$, $a_{1}$, \ldots, $a_{p-1}$, fixed points
of $F_{0}$, $F_{1}$, \ldots, $F_{p-1}$, respectively, that is
\begin{align*}
F_{0}\left(  a_{0}\right)   &  =a_{0},\\
F_{1}\left(  a_{1}\right)   &  =a_{1},\\
&  \vdots\\
F_{j}\left(  a_{j}\right)   &  =a_{j},\\
&  \vdots\\
F_{p-1}\left(  a_{p-1}\right)   &  =a_{p-1}.
\end{align*}

\item \label{DD2}The non-hyperbolicity condition%
\[
\left.  \partial_{x}F_{j}\left(  x\right)  \right\vert _{x=a_{j}}%
=\underset{i=0}{\overset{p-1}{%
%TCIMACRO{\dprod }%
%BeginExpansion
{\displaystyle\prod}
%EndExpansion
}}\partial_{x}f_{i}\left(  a_{j}\right)  =1.
\]

\item \label{DD3}Higher degeneracy conditions%
\[
\ \left.  \partial_{x^{i}}F_{j}\left(  x\right)  \right\vert _{x=a_{j}%
}=0,\mathit{\ }\text{\textit{for every }}2\leq i\leq\mu.
\]

\item \label{DD4}The non-degeneracy condition%
\[
\ \left.  \partial_{x^{\mu+1}}F_{j}\left(  x\right)  \right\vert _{x=a_{j}%
}\not =0.
\]

\end{enumerate}

Then, all the compositions $F_{m}$, $0\leq m\leq p-1$, satisfy%
\[
\left.  \partial_{x^{i}}F_{m}\left(  x\right)  \right\vert _{x=a_{m}}=0,\text{
for every }2\leq i\leq\mu
\]
and
\[
\left.  \partial_{x^{\mu+1}}F_{m}\left(  x\right)  \right\vert _{x=a_{m}%
}\not =0.
\]

\end{thm}

\begin{pf}
Without loss of generality we consider that the hypothesis apply when $j=0$,
what can be done re-indexing the maps of the $p$-periodic system. We now apply
Theorem \ref{t1} to the alternating system $f_{0}$, $f_{p-1}\circ\ldots\circ
f_{1}$ with compositions $F=F_{0}$ and $G=F_{1}=f_{0}\circ F\circ f_{0}^{-1}$,
making $a=a_{0}$, $b=a_{1}$ and getting the result for $j=1$. Applying the
same argument repeatedly, the result follows immediately for all the cyclic
compositions $F_{m}$, $0\leq m\leq p-1$.
\end{pf}

\begin{rmk}
\label{rem}The same result holds for $k$-periodic points of the compositions
$F_{j}$, i.e., $k\times p$-periodic points of the alternating system, since in
that case we apply Theorem \ref{t1} to the alternating system with
compositions $F=F_{0}^{k}\ $and $G=F_{1}^{k}=f_{0}\circ F^{k}\circ f_{0}^{-1}$.
\end{rmk}

After this result we can choose the composition order that makes the
bifurcation equations easier to solve.

\subsection{Invariance of the transversality
conditions\label{unfoldingsection}}

\subsubsection{Alternating maps}

In this paragraph we prove the symmetry for the transversality conditions
concerning cyclic compositions of $p$ maps. We follow the same technique of
proving the result for the alternating maps\footnote{Notation that we adopt in
this paragraph simplify the presentation of the next results and proofs, we
replace $f_{0}$ by $f$ and $f_{1}$ by $g$. The compositions are $F=g\circ f$
and $G=f\circ g$.} $f$ and $g$ and generalizing it to periodic points of
$p$-periodic systems. Suppose that there exists the solution $\Lambda_{0}$ in
the parameter space of the bifurcation equations \ref{bif} and that the
non-degeneracy condition \ref{nond} holds at $\Lambda_{0}$, that solution
coexists with fixed points $a$ and $b$ for the compositions $F$ and $G$.

Consider the map $\mathcal{F}=\left(  F,F_{x},F_{xx},\ldots,F_{x^{\mu-1}%
}\right)  $ with the derivatives of the composition $F$ and the Jacobian
determinant of $\mathcal{F}$, now as a function of $\Lambda$. It is the
determinant%
\begin{equation}
J_{\Lambda}\mathcal{F}\left(  x,\Lambda\right)  =\det\left[
\begin{array}
[c]{lll}%
F_{\lambda_{1}}(x,\Lambda) & \cdots & F_{\lambda_{\mu}}(x,\Lambda)\\
\vdots & \ddots & \vdots\\
F_{x^{\mu-1}\lambda_{1}}(x,\Lambda) & \cdots & F_{x^{\mu-1}\lambda_{\mu}%
}(x,\Lambda)
\end{array}
\right]  . \label{Wronsk}%
\end{equation}
Consider similar definitions for $\mathcal{G}\left(  x,\Lambda\right)  $ and
$J_{\Lambda}\mathcal{G}\left(  x\right)  $ relative to the composition $G$.
With the previous definitions we can establish the next lemma.

\begin{lem}
\label{Nondegenaracy}The Jacobians $J_{\Lambda}\mathcal{F}\left(
x,\Lambda\right)  $ and $J_{\Lambda}\mathcal{G}\left(  x,\Lambda\right)  $
computed at the solutions of the bifurcation conditions (\ref{bif}) and
(\ref{nond}) satisfy the equality%
\begin{equation}
J_{\Lambda}\mathcal{F}\left(  a,\Lambda_{0}\right)  =\left(  f_{x}\left(
a,\Lambda_{0}\right)  \right)  ^{^{\frac{3\mu-\mu^{2}}{2}}}J_{\Lambda
}\mathcal{G}\left(  b,\Lambda_{0}\right)  . \label{Unfold}%
\end{equation}

\end{lem}

\begin{pf}
The proof rests on the fact that we can obtain the lines of the Jacobian
matrix $\left[  J_{\Lambda}\mathcal{F}\left(  a,\Lambda_{0}\right)  \right]  $
of $\mathcal{F}$ using Gaussian manipulation of the Jacobian matrix $\left[
J_{\Lambda}\mathcal{G}\left(  b,\Lambda_{0}\right)  \right]  $ of
$\mathcal{G}$.

Consider%
\[
\left[  J_{\Lambda}\mathcal{F}\left(  a,\Lambda_{0}\right)  \right]  =\left[
\begin{array}
[c]{c}%
L_{1}\\
L_{2}\\
\vdots\\
L_{\mu}%
\end{array}
\right]  ,
\]
where $L_{i}$ denotes the $i$ line of the matrix. We have to prove that
\[
\left[  J_{\Lambda}\mathcal{G}\left(  b,\Lambda_{0}\right)  \right]  =\left[
\begin{array}
[c]{c}%
\alpha_{11}L_{1}\\
\alpha_{21}L_{1}+\alpha_{22}L_{2}\\
\vdots\\%
%TCIMACRO{\dsum \limits_{j=1}^{\mu}}%
%BeginExpansion
{\displaystyle\sum\limits_{j=1}^{\mu}}
%EndExpansion
\alpha_{\mu j}L_{j}%
\end{array}
\right]  ,
\]
i.e.%
\[
\left[  J_{\Lambda}\mathcal{G}\left(  b,\Lambda_{0}\right)  \right]  =\left[
\begin{array}
[c]{cccc}%
\alpha_{11} & 0 & \cdots & 0\\
\alpha_{21} & \alpha_{22} & \cdots & 0\\
\vdots & \vdots & \ddots & \vdots\\
\alpha_{\mu1} & \alpha_{\mu2} & \cdots & \alpha_{\mu\mu}%
\end{array}
\right]  \left[
\begin{array}
[c]{c}%
L_{1}\\
L_{2}\\
\vdots\\
L_{\mu}%
\end{array}
\right]  ,
\]
which is
\[
\left[  J_{\Lambda}\mathcal{G}\left(  b,\Lambda_{0}\right)  \right]  =A\left[
J_{\Lambda}\mathcal{F}\left(  a,\Lambda_{0}\right)  \right]  \text{,}%
\]
where
\[
A=\left[
\begin{array}
[c]{cccc}%
\alpha_{11} & 0 & \cdots & 0\\
\alpha_{21} & \alpha_{22} & \cdots & 0\\
\vdots & \vdots & \ddots & \vdots\\
\alpha_{\mu1} & \alpha_{\mu2} & \cdots & \alpha_{\mu\mu}%
\end{array}
\right]
\]
and the central point that $\det A$ must be different from $0$.

The fact that for general $\mu\geq1$ the matrix $A$ is a lower triangular
matrix is trivial. Bellow, we prove that each entry of the main diagonal is%
\begin{equation}
\alpha_{jj}=\left(  f_{x}\left(  a,\Lambda_{0}\right)  \right)  ^{2-j}\text{,
}1\leq j\leq\mu\text{,} \label{diagonal}%
\end{equation}
this equality implies that all such entries are different from zero after the
non-hyperbolicity condition at the bifurcation%
\[
\left.  \partial_{x}F\left(  x,\Lambda\right)  \right\vert _{x=a,\Lambda
=\Lambda_{0}}=g_{x}\left(  b,\Lambda_{0}\right)  f_{x}\left(  a,\Lambda
_{0}\right)  =1\text{.}%
\]
Moreover, (\ref{diagonal}) implies that the determinant of $A$ is%
\[
\det A=%
%TCIMACRO{\dprod \limits_{j=1}^{\mu}}%
%BeginExpansion
{\displaystyle\prod\limits_{j=1}^{\mu}}
%EndExpansion
\left(  f_{x}\left(  a,\Lambda_{0}\right)  \right)  ^{2-j}=\left(
f_{x}\left(  a,\Lambda_{0}\right)  \right)  ^{^{\frac{3\mu-\mu^{2}}{2}}%
}\text{,}%
\]
as desired.

We prove now equality (\ref{diagonal}). Note that $G\circ f=f\circ g\circ
f=f\circ F$. We derive this local conjugacy in order to $\lambda_{i}$, with
$i=1,2,\ldots,\mu$. We have
\begin{equation}
\partial_{\lambda_{i}}G\left(  f\left(  x,\Lambda\right)  ,\Lambda\right)
=\partial_{\lambda_{i}}f\left(  F\left(  x,\Lambda\right)  ,\Lambda\right)
\text{,} \label{Original}%
\end{equation}
with%
\begin{equation}%
\begin{array}
[c]{l}%
\partial_{\lambda_{i}}G\left(  f\left(  x,\Lambda\right)  ,\Lambda\right)
=G_{\lambda_{i}}\left(  f\left(  x,\Lambda\right)  ,\Lambda\right)
+G_{x}\left(  x,\Lambda\right)  f_{\lambda_{i}}\left(  x,\Lambda\right)
\text{,}\\
\partial_{\lambda_{i}}f\left(  F\left(  x,\Lambda\right)  ,\Lambda\right)
=f_{\lambda_{i}}\left(  x,\Lambda\right)  +f_{x}\left(  F\left(
x,\Lambda\right)  ,\Lambda\right)  F_{\lambda_{i}}\left(  x,\Lambda\right)
\text{.}%
\end{array}
\label{conjderiv}%
\end{equation}
At the points $\left(  a,\Lambda_{0}\right)  $ and $\left(  b,\Lambda
_{0}\right)  $ equating the second members of (\ref{conjderiv})\ one has%
\[
G_{\lambda_{i}}\left(  b,\Lambda_{0}\right)  +G_{x}\left(  b,\Lambda
_{0}\right)  f_{\lambda_{i}}\left(  a,\Lambda_{0}\right)  =f_{\lambda_{i}%
}\left(  a,\Lambda_{0}\right)  +f_{x}\left(  a,\Lambda_{0}\right)
F_{\lambda_{i}}\left(  a,\Lambda_{0}\right)  \text{,}%
\]
which after conditions (\ref{bif}) is at the bifurcation point%
\begin{equation}
G_{\lambda_{i}}\left(  b,\Lambda_{0}\right)  =f_{x}\left(  a,\Lambda
_{0}\right)  F_{\lambda_{i}}\left(  a,\Lambda_{0}\right)  \text{,}
\label{First Line}%
\end{equation}
this equality gives the relation between the first rows of the Jacobians.

To get the relations between the second rows we consider the derivative
relative to $x$ of (\ref{conjderiv}), we present only the terms that matter
for the computation of the main diagonal of $A$
\begin{equation}%
\begin{array}
[c]{l}%
\partial_{\lambda_{i}x}G\left(  f\left(  x,\Lambda\right)  ,\Lambda\right)
=G_{\lambda_{i}x}\left(  f\left(  x,\Lambda\right)  ,\Lambda\right)
f_{x}\left(  x,\Lambda\right)  +\cdots\\
\\
\partial_{\lambda_{i}x}f\left(  F\left(  x,\Lambda\right)  ,\Lambda\right)
=\cdots+f_{x}\left(  F\left(  x,\Lambda\right)  ,\Lambda\right)
F_{\lambda_{i}x}\left(  x,\Lambda\right)  \text{.}%
\end{array}
\label{conjdreiv1}%
\end{equation}
which after conditions (\ref{bif}) gives at the bifurcation value equalizing
the right hand sides of (\ref{conjdreiv1})
\begin{equation}
G_{x\lambda_{i}}\left(  b,\Lambda_{0}\right)  =l.o.t+F_{x\lambda_{i}}\left(
a,\Lambda_{0}\right)  \text{,} \label{Second Line}%
\end{equation}
where $l.o.t$ stands for \textquotedblleft lower order terms\textquotedblright%
\ in terms of derivatives on the dynamical variable of $G_{\lambda_{i}}$ and
$F_{\lambda_{i}}$, terms that do not appear in the main diagonal of $A$. This
expression gives the relation between the second rows of the Jacobians.

To get the relation between the third rows of the two Jacobians, we consider
the derivative of (\ref{conjdreiv1}) regarding $x$%
\[%
\begin{array}
[c]{l}%
\partial_{\lambda_{i}x^{2}}G\left(  f\left(  x,\Lambda\right)  ,\Lambda
\right)  =G_{\lambda_{i}x^{2}}\left(  f\left(  x,\Lambda\right)
,\Lambda\right)  f_{x}^{2}\left(  x,\Lambda\right)  +\cdots\\
\\
\partial_{\lambda_{i}x^{2}}f\left(  F\left(  x,\Lambda\right)  ,\Lambda
\right)  =\cdots+f_{x}\left(  F\left(  x,\Lambda\right)  ,\Lambda\right)
F_{\lambda_{i}x^{2}}\left(  x,\Lambda\right)  \text{,}%
\end{array}
\]
which after conditions(\ref{bif}) gives at the bifurcation value%
\[
G_{x^{2}\lambda_{i}}\left(  b,\Lambda_{0}\right)  f_{x}^{2}\left(
a,\Lambda_{0}\right)  =l.o.t+f_{x}\left(  a,\Lambda_{0}\right)  F_{x^{2}%
\lambda_{i}}\left(  a,\Lambda_{0}\right)  \text{,}%
\]
repeating this process and using the Fa\`{a} di Bruno Formula (\ref{bruno})
and the bifurcation equations (\ref{bif}), knowing that the lower order terms
in derivatives relative to $x$ (order less than $j-1$) do not contribute to
the diagonal of $A$, we have for $1\leq j\leq\mu$%
\[
G_{x^{j-1}\lambda_{i}}\left(  b,\Lambda_{0}\right)  \left(  f_{x}\left(
a,\Lambda_{0}\right)  \right)  ^{j-1}=l.o.t+f_{x}\left(  a,\Lambda_{0}\right)
F_{x^{j-1}\lambda_{i}}\left(  a,\Lambda_{0}\right)  \text{,}%
\]
dividing by $\left(  f_{x}\left(  a,\Lambda_{0}\right)  \right)  ^{j-1}$,
which can not be zero from the second line of equations (\ref{bif}), we
obtain
\[
G_{x^{j-1}\lambda_{i}}\left(  b,\Lambda_{0}\right)  =l.o.t+\left(
f_{x}\left(  a,\Lambda_{0}\right)  \right)  ^{2-j}F_{x^{j-1}\lambda_{i}%
}\left(  a,\Lambda_{0}\right)  \text{,}%
\]
the desired result.
\end{pf}

\subsubsection{Cyclic composition of $p$ maps}

In the general setting of $p$-periodic maps and considering the last lemma,
the transversality conditions for $A_{\mu}$ bifurcations of $k\times p$
periodic points of the first two possible compositions of maps are such that%
\begin{equation}
J_{\Lambda}\mathcal{F}_{0}^{k}\left(  a_{0},\Lambda_{0}\right)  \not =%
0\Rightarrow J_{\Lambda}\mathcal{F}_{1}^{k}\left(  a_{1},\Lambda_{0}\right)
\not =0, \label{Unfold2}%
\end{equation}
where $J_{\Lambda}$ was defined in (\ref{Wronsk}), because $F_{0}^{k}$ and
$F_{1}^{k}$ are the two compositions of alternating maps as we have seen in
Remark \ref{rem}. The generalization to periodic points of all the cyclic
compositions of $p$-periodic maps poses no difficulties and the proof is
obtained by repeated application of Lemma \ref{Nondegenaracy}.

\begin{thm}
\label{Unfo}If one of the transversality conditions of the $A_{\mu}$
bifurcation for $k\times p$-periodic orbits of $p$-periodic maps at
$\Lambda=\Lambda_{0}$ is satisfied, say%
\[
J_{\Lambda}\mathcal{F}_{j}^{k}\left(  a_{j},\Lambda_{0}\right)  \not =0
\]
then it is satisfied for all the cyclic compositions of the individual maps,
i.e.,
\begin{align*}
J_{\Lambda}\mathcal{F}_{0}^{k}\left(  a_{0},\Lambda_{0}\right)   &
\not =0\text{,}\\
\vdots & \\
J_{\Lambda}\mathcal{F}_{p-1}^{k}\left(  a_{p-1},\Lambda_{0}\right)   &
\not =0\text{.}%
\end{align*}

\end{thm}

\subsection{Conclusion}

The invariance of degeneracy and transversality conditions, imply that the
bifurcation problem with $\mu$ degeneracy conditions on the iteration variable
is independent on the choice of the cyclic order in the composition of the
maps when the maps are sufficiently differentiable. This invariance is
fundamental, since it means that we can define local bifurcations in a unique
way for families of $p$-periodic maps using any of the compositions of the
particular maps.

In particular, the bifurcation set in the parameter space is the same for all
$F_{j}^{k}$.

The main conclusion of this study, is that it suffices to solve the
bifurcation conditions applied to one of the $F_{j}^{k}$ possible compositions
to obtain the bifurcation set. The bifurcation is well defined using the
bifurcation conditions on the composition families. Each fiber replicates the
behavior of the others. Hence, the local bifurcations studied in this work of
$p$-periodic difference equations are defined by the same rules of the local
bifurcations of autonomous systems.

\section{\label{Examples}Examples}

We conclude this work with the study of the particular case of alternating
maps. We establish some useful criteria about the existence of the swallowtail
singularity for alternating maps and give two examples exhibition this type of bifurcation.

Let $f\in{\mathcal{C}}^{3}\left(  I\right)  $. The Schwarzian derivative of
$f$ is%
\[
Sf(x)=\frac{f^{\prime\prime\prime}\left(  x\right)  }{f^{\prime}\left(
x\right)  }-\frac{3}{2}\left(  \frac{f^{\prime\prime}\left(  x\right)
}{f^{\prime}\left(  x\right)  }\right)  ^{2},
\]
defined for every $x$ in $I$, that is not a critical point of $f$.

\begin{prop}
\label{SchwarzSign}Consider\ an alternating system with families, $f=f_{0}$
and $g=f_{1}$, satisfying all the conditions of the $A_{3}$ bifurcation, i.e.,
swallowtail bifurcation, together with the transversality conditions. If one
of the maps say, without loss of generality, $g$, has Schwarzian derivative
different from zero at $b$, $Sg\left(  b\right)  \not =0$, then the product of
the Schwarzian derivatives must be negative at the swallowtail bifurcation
point, i.e.,%
\[
Sg\left(  b\right)  \cdot Sf\left(  a\right)  <0\text{.}%
\]

\end{prop}

\begin{pf}
Recall the example \ref{example} with the same notation for derivatives.
Consider $f$ and $g$ in the conditions of that example. Remember the
equalities obtained and the same simplifying notation. If $g$ is assumed to
have negative Schwarzian derivative at $b$, one has%
\[
\frac{g_{3}}{g_{1}}-\frac{3}{2}\left(  \frac{g_{2}}{g_{1}}\right)  ^{2}<0,
\]
then, from (\ref{cc2}) and (\ref{g3}) one has%
\[
\frac{g_{3}}{g_{1}}=\frac{3\left(  f_{1}\right)  ^{2}}{\left(  f_{1}\right)
^{4}}-\frac{f_{3}}{\left(  f_{1}\right)  ^{3}}=\frac{1}{\left(  f_{1}\right)
^{2}}\left(  3\left(  \frac{f_{2}}{f_{1}}\right)  ^{2}-\frac{f_{3}}{f_{1}%
}\right)  <\frac{3}{2}\left(  \frac{g_{2}}{g_{1}}\right)  ^{2},
\]
and, by the equality (\ref{g2})
\[
\frac{g_{2}}{g_{1}}=-\frac{f_{2}}{\left(  f_{1}\right)  ^{2}},
\]
the inequality above becomes%
\[
3\left(  \frac{f_{2}}{f_{1}}\right)  ^{2}-\frac{f_{3}}{f_{1}}<\frac{3}%
{2}\left(  \frac{f_{2}}{f_{1}}\right)  ^{2}.
\]
Therefore,
\[
\frac{f_{3}}{f_{1}}-\frac{3}{2}\left(  \frac{f_{2}}{f_{1}}\right)  ^{2}>0.
\]
Therefore, if $g$ has negative Schwarzian derivative at $b$, then $f$ must
have positive Schwarzian derivative at $a$.

On the other hand, if $f$ is assumed to have negative Schwarzian derivative,
then by similar reasonings
\[
\frac{f_{3}}{f_{1}}=-\frac{g_{3}}{g_{1}}\left(  f_{1}\right)  ^{2}+3\left(
\frac{f_{2}}{f_{1}}\right)  ^{2}<\frac{3}{2}\left(  \frac{f_{2}}{f_{1}%
}\right)  ^{2}.
\]
This implies%
\[
\frac{g_{3}}{g_{1}}\left(  f_{1}\right)  ^{2}>\frac{3}{2}\left(  \frac{f_{2}%
}{f_{1}}\right)  ^{2}%
\]
and, as%
\[
\frac{g_{2}}{g_{1}}=-\frac{f_{2}}{\left(  f_{1}\right)  ^{2}},
\]
it follows that
\[
\frac{g_{3}}{g_{1}}>\frac{3}{2}\left(  \frac{f_{2}}{\left(  f_{1}\right)
^{2}}\right)  ^{2}=\frac{3}{2}\left(  \frac{g_{2}}{g_{1}}\right)  ^{2}.
\]
Hence, if $f$ has negative Schwarzian derivative then $g$ must have positive
Schwarzian derivative.
\end{pf}

As in \cite{EH}, while working on pitchfork bifurcation, we can state the two
following propositions for the $A_{3}$ degenerate bifurcation. The proofs are similar.

\begin{prop}
Let $f$ and $g$ be ${\mathcal{C}}^{3}$ alternating maps. If $f$ is strictly
increasing and $g$ is strictly decreasing in $x$ (analogously, if $f$ is
strictly decreasing and $g$ is strictly increasing), then the alternating
system associated with $f$ and $g$ cannot have a $A_{3}$ degenerate bifurcation.
\end{prop}

\begin{prop}
Let $f$ and $g$ be ${\mathcal{C}}^{3}$ alternating maps. If $f$ and $g$ are
both strictly increasing in the dynamic variable and one of these two
situations happens
\[
\min_{x\in I_{0},I_{1}}(\partial_{x^{2}}f(x),\partial_{x^{2}}%
g(x))>0\ \text{(both\/convex)}%
\]
or
\[
\max_{x\in I_{0},I_{1}}(\partial_{x^{2}}f(x),\partial_{x^{2}}%
g(x))<0\ \text{(both\/concave),}%
\]
then the alternating system generated by them cannot have a swallowtail bifurcation.
\end{prop}

\begin{prop}
If $f$ and $g$ are both strictly decreasing in the dynamic variable and this
situation happens
\[
\max_{x\in I_{0},I_{1}}(\partial_{x^{2}}f(x),\partial_{x^{2}}%
g(x))>0\text{\ and \/}\min_{x\in I_{0},I_{1}}(\partial_{x^{2}}f(x),\partial
_{x^{2}}g(x))<0,
\]
i.e., one is concave and one is convex or vice-versa, then the alternating
system generated by them cannot have a $A_{3}$ degenerate bifurcation.
\end{prop}

At this point, we consider two concrete examples. The first one is an
alternating system with polynomial families, and the second with the tangent
family, $\lambda\tan x$, and a polynomial family. The first example is
relatively easy to compute, but the second one has elusive roots due to its
high degeneracy. Thus, some numeric work was necessary. We present only the
solutions and discard the tedious computations of the second example.

\begin{figure}[h]
\begin{center}
\includegraphics[
height=2.5755in,
width=2.5755in
]{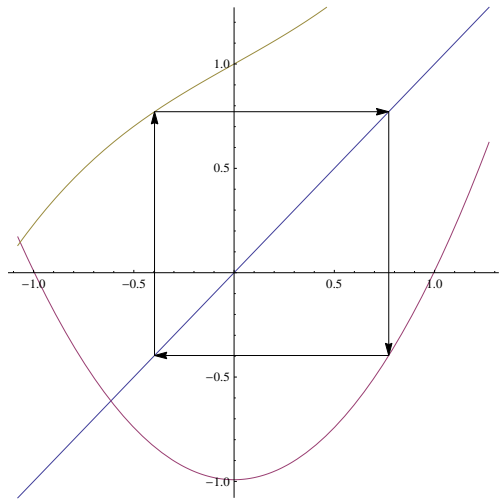}
\end{center}
\caption{The geometry of the individual maps at the Swallowtail bifurcation
point $A_{3}$ for example \ref{Example1}. Please note that one map is convex
and the other is concave.}%
\label{FIG8}%
\end{figure}

\begin{exmp}
\label{Example1}\label{Example1 copy(1)}Alternating system $f$ and $g$ with a
quadratic polynomial $f_{0}=x^{2}+\lambda_{1}$ and a cubic polynomial
$f_{1}=\lambda_{3}x^{3}+\lambda_{2}x+1$, defined in the real line. The
compositions are%
\[
F\left(  x\right)  =f_{1}\circ f_{0}\left(  x\right)  =\lambda_{3}%
x^{6}+3\lambda_{3}\lambda_{1}x^{4}+\left(  3\lambda_{3}\lambda_{1}^{2}%
+\lambda_{2}\right)  x^{2}+1+\lambda_{2}\lambda_{1}+\lambda_{3}\lambda_{1}^{3}%
\]
and
\[
G\left(  x\right)  =f_{0}\circ f_{1}\left(  x\right)  =\lambda_{3}^{2}%
x^{6}+2\lambda_{2}\lambda_{3}x^{4}+2\lambda_{3}x^{3}+\lambda_{2}^{2}%
x^{2}+2\lambda_{2}x+1+\lambda_{1}\text{.}%
\]
The bifurcation equations (\ref{bif}) for $F$ or $G$ have solutions%
\[
\text{ }\lambda_{1}=-\frac{3^{5}}{5\cdot7^{2}}\text{, }\lambda_{2}=\frac
{5^{2}\cdot7}{2^{2}\cdot3^{4}}\text{ and }\lambda_{3}=\frac{7^{5}}{2^{4}%
\cdot3^{10}}\text{,}%
\]
with
\[
a=\frac{3^{3}}{5\cdot7}\text{, }b=-\frac{2\cdot3^{5}}{5^{2}\cdot7^{2}}\text{,}%
\]
such that $f_{0}\left(  a\right)  =b$ and $f_{1}\left(  b\right)  =a$. The
Schwarzian derivatives are%
\[
Sf\left(  a\right)  =\frac{1}{b}=-\frac{5^{2}\cdot7^{2}}{2\cdot3^{5}}\text{,
}Sg\left(  b\right)  =\frac{1}{6}\frac{1}{b^{2}}=\frac{5^{4}\cdot7^{4}}%
{2^{3}\cdot3^{11}}\text{,}%
\]
naturally, with opposite signs at the bifurcation points, accordingly with
Proposition \ref{SchwarzSign}. Obviously, $SF\left(  a\right)  =SG\left(
b\right)  =0$ at the bifurcation points.

The transversality condition (\ref{Unfold2}) is from Lemma \ref{Nondegenaracy}
equal for $F$ and $G$ (since $\frac{3.3-3^{2}}{2}=0$) and gives
\[
J_{\Lambda}\mathcal{F}\left(  a,\Lambda_{0}\right)  =J_{\Lambda}%
\mathcal{G}\left(  b,\Lambda_{0}\right)  =-\frac{2^{4}\cdot3^{10}}{5^{4}%
\cdot7^{3}}\not =0\text{.}%
\]

We can see at Figure \ref{FIG8} the geometry of the individual maps at the
$A_{3}$ singularity. Both functions are increasing at suitable neighborhoods
of $a$ and $b$ and one function is concave and the other is convex.

The bifurcation set is exactly similar to the one depicted in Figure
\ref{FIG9}. We have the same behavior of the two possible compositions $F$ and
$G$.
\end{exmp}

\begin{figure}[h]
\begin{center}
\includegraphics[
height=2.5639in,
width=2.6941in
]{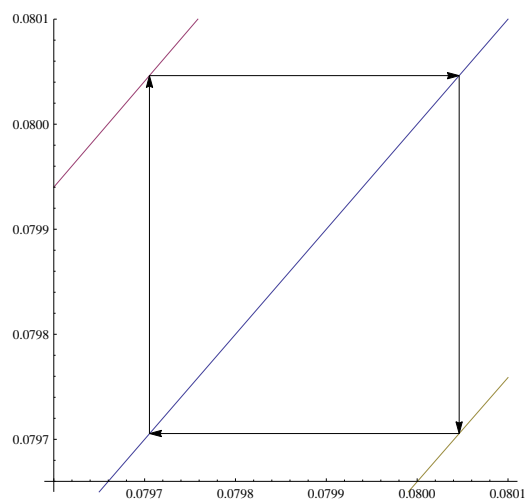}
\end{center}
\caption{The orbit of the bifurcation points $a$ and $b$ of \ $f_{0},f_{1}$ in
example \ref{EX2} viewed as a cobweb diagram. The maps, one concave and the
other convex, are almost parallel.}%
\label{FIG12}%
\end{figure}

\begin{exmp}
\label{EX2}Consider now the family of real alternating maps $f_{0}$ and
$f_{1}$ with $f_{0}\left(  x\right)  =-x^{4}+\lambda_{1}x^{2}+x+\lambda_{2}$
and $f_{1}\left(  x\right)  =\lambda_{3}\tan x$, defined in suitable open sets
near the solutions of the swallowtail bifurcation equations. We have the
solutions of the bifurcation conditions $a\simeq0.0797053$, $b\simeq
0.0793675$, $\lambda_{1}\simeq-0.0400839$, $\lambda_{2}\simeq-0.0000428492$,
$\lambda_{3}\simeq1.00215$. The non-degeneracy condition gives $F_{x^{4}%
}\left(  a\right)  \simeq-26.7$. We note that $a$ and $b$ are very near each
other and the maps are almost parallel at the bifurcation points, as we can
see at Figure \ref{FIG12}.The Schwarzian derivatives are
\[
Sf\left(  a\right)  =-1.96648,\text{ }Sg\left(  b\right)  =2.
\]
The transversality condition (\ref{Unfold2}) is again equal for $F$ and $G$
and is
\[
J_{\Lambda}\mathcal{F}\left(  a,\Lambda_{0}\right)  =J_{\Lambda}%
\mathcal{G}\left(  b,\Lambda_{0}\right)  =-2.08013\text{.}%
\]

\end{exmp}

These two examples, exhibiting the swallowtail bifurcation, produce evidence
that high degeneracy bifurcations can occur in concrete examples and the
theory is not void of applications.

\textbf{Acknowledgement} The author wants to thank the anonymous referee about
his precious remarks and suggestions that improved a great deal the revised
version of this paper. The author also thanks Michal Misiurewicz for the
fruitful discussion in the Polish countryside about the proof of Lemma
\ref{lema}. The author was partially funded through project
PEst-OE/EEI/LA0009/2013 for CAMGSD.

\bibliographystyle{elsart-num-sort}
\bibliography{BibloH2}

\end{document}